\documentclass[12pt]{article}
\usepackage{amsmath,amsfonts,amssymb,amsthm,amscd}
\title{Forking in the free group}
\author{Anand
Pillay\thanks{Supported by a Marie Curie Chair MEXC -CT-2005 -024052}\\School of Mathematics\\University of 
Leeds\\Leeds LS2 9JT\\UK\\pillay@maths.leeds.ac.uk}
\newtheorem{Theorem}{Theorem}[section]
\newtheorem{Proposition}[Theorem]{Proposition}
\newtheorem{Definition}[Theorem]{Definition}
\newtheorem{Remark}[Theorem]{Remark}
\newtheorem{Lemma}[Theorem]{Lemma}
\newtheorem{Corollary}[Theorem]{Corollary}
\newtheorem{Fact}[Theorem]{Fact}

\newtheorem{Problem}[Theorem]{Problem}
\newcommand{\R}{\mathbb R}
\newcommand{\Q}{\mathbb Q}
\newcommand{\Z}{\mathbb Z}
\begin{document}
\maketitle

\begin{abstract}
We study model-theoretic and stability-theoretic properties of the nonabelian free group in the light of Sela's recent result
\cite{Sela VIII} on stability and results announced by Bestvina and Feighn on ``negligible subsets" of free groups. 
We point out analogies
between the free group and so-called bad groups of finite
Morley rank, and prove ``non $CM$-triviality" of the free
group.
\end{abstract}

\section{Introduction}
Let $F_{n}$ denote the free group on $n$ generators. We view
$F_{n}$ as a first order structure $(F_{n},\cdot,^{-1},1)$
in the language $L$ of groups, where $\cdot$ is the group
operation, $^{-1}$ inversion, and $1$ is the identity
element. In a recent preprint \cite{Sela VIII} Zlil Sela proves the rather astounding result:
\newline
{\em Theorem (A).}  For any $n$, $Th(F_{n})$ is stable.

\vspace{2mm}
\noindent
This built on  a sequence of papers culminating in \cite{Sela V-1}, \cite{Sela V-2}, \cite{Sela VI},  which included the results:
\newline
{\em Theorem (B).} (i) for any $2\leq m < n$, the natural embedding of $F_{m}$ in $F_{n}$ is an elementary embedding, and
\newline
(ii) The common (complete) theory $T_{fg}$ of the groups $F_{n}$ ($n\geq 2$) has quantifier elimination down to Boolean combinations of $\forall\exists$ formulas.

\vspace{2mm}
\noindent
Other stable groups (that is groups whose first order theory
is stable) are (i) any commutative group (in the group
language), and (ii) the group $G(K)$ of $K$-points of an
algebraic group $G$ over an algebraically closed field $K$,
where now the language is the Zariski-language, namely we
have relations or predicates for all Zariski-closed subsets
of $G^{n}(K)$ all $n$.

There is in place a beautiful theory of stable groups and group
actions, ``equivariant stability theory", which is due in
full generality to Bruno Poizat \cite{Poizat} and \cite{Poizat-book}, (building on earlier work of Macintyre, Zilber, Cherlin and Shelah among others) and which consciously
borrows terminology such as stabilizers, 
connected components, and generics from the theory of algebraic groups.
In particular this theory yields a
well-behaved notion of genericity or largeness for {\em definable}
subsets of a stable group. On the other hand, as we heard
first from Sela in 2003, Bestvina and Feighn have come up
with a specific combinatorial notion of largeness for
arbitrary subsets of a free group. We point out in section 2
that (unpublished) results announced by Bestvina and Feighn
imply that these two  notions of largeness, stability-theoretic and
combinatorial, coincide for {\em
definable} subsets of the free group. We take the
opportunity in the rest of section 2 to point out several other consequences, 
some of
which are already known and were even pointed out to me by Sela, such as the structure of definable
subgroups.

The main point of the paper is in section 3 where we
show that the free group is rather complicated from the
point of view of the ``geometry of forking".  This notion of
complicatedness ``non $CM$-triviality" or ``$2$-ampleness"
is rather delicate, and was originally defined in
\cite{Hrushovski} to describe properties of dimension (rank)
and algebraic closure in certain {\em strongly minimal}  or
finite Morley rank structures. Roughly speaking
$CM$-triviality forbids the existence of a certain definable
``point-line-plane" configuration, much as $1$-basedness
forbids a certain definable ``point-line" configuration.
Hrushovski \cite{Hrushovski} gave a counterexample to a conjecture of Zilber
by constructing a strongly minimal set which is not
$1$-based but does not interpret an infinite field. He also
observed that his new strongly minimal set is in fact
$CM$-trivial.  It remains an important open question whether
there is a non $CM$-trivial strongly minimal theory which
does not interpret an infinite field.

The free group is very far from having finite Morley rank.
It is not even superstable. But the notion of
$CM$-triviality still makes sense (as does
$1$-basedness) in arbitrary stable structures.  Hrushovski's
general method of construction, via $\delta$-functions,
strong embeddings and amalgation, can be used to produce new
stable structures, not necessarily of finite rank, and again
these are typically $CM$-trivial structures (or
$CM$-trivial over the starting data). So the heuristic
conclusion is that the free group cannot arise from such a
Hrushovski construction. On the other hand it is not so very
hard to produce non $CM$-trivial stable structures (even
$\omega$-stable ones) which do not interpret infinite
fields. One such structure, the free pseudospace
  was constructed by rather ad hoc
means by Baudisch and the author \cite{B-P}. It turned out to
be closely related to an unpublished example of a {\em
non-equational} $
\omega$-stable theory due to Hrushovski and Srour \cite{Hrushovski-Srour}.  Other
examples were constructed in a rather more systematic fashion
by David Evans \cite{Evans}.  Anyway the fact that such an
example (namely the free group) occurs in nature is rather
interesting and the reason for writing this paper. 

In \cite{Pillay1} we proved that any simple noncommutative
group of finite Morley rank is non $CM$-trivial. A key case
to deal with was that of so-called bad groups (simple groups
in which Borels are nilpotent).  Since learning of Sela's
work on the stability of the free group we wanted to carry
over this proof of non $CM$-triviality to the free group,
but there were certain technical obstructions due to being
outside the finite Morley rank context. The results
announced by Bestvina and Feighn  give an alternative
computational tool, and this is what we explain in section 3.

We conclude in section 4 with  some natural questions about
the model theory of the free group.

In this paper we will freely use the language, notions, and techniques of model theory and stability theory. References are \cite{Marker}, \cite{Poizat-book}, as well as \cite{Pillay-book}. However for the benefit of the more general reader we will take the opportunity in the rest of this introduction to explain something about {\em stable groups}.

As a matter of notation, by a definable set in a first order structure $M$, we mean a subset of some $M\times..\times M$ which is definable possibly with parameters in $M$. If we want to specify that the defining parameters come from a subset $A$ of $M$ we say $A$-definable. So $\emptyset$-definable means definable without parameters.

By a group in the sense of model theory, we usually mean a group $(G,\cdot)$ equipped possibly with additional relations or preducates, namely some subsets $R_{i}$ of cartesian powers $G^{n_{i}}$ for $i$ ranging over an index set $I$.
For example if $G$ is an algebraic group over an algebraically closed field $k$, it is natural to equip $G(k)$ with predicates for all Zariski closed subsets of its Cartesian powers.  If $G$ is a semialgebraic real Lie group (like the connected component of $GL(n,\R)$), it would be natural to equip $G$ with all semialgebraic subsets of its Cartesian powers. But even if we consider $G$ without explicit additional relations (as we do the free group), all the subsets of $G,G\times G,...$ {\em first order definable} in $(G,\cdot)$, will be part of the structure and may be rather complicated (as in the case of the free group).
In any case the point here is that model-theorists often treat groups as if they were objects  of geometry, like algebraic groups and Lie groups. 

Given a complete theory $T$, we are interested in definable sets in arbitrary models of $T$, in particular in {\em saturated} models of $T$. Likewise for the free groups $F$ we are typically interested not only in $F$ but in elementary extensions of $F$.

By a {\em stable group} we mean a group $(G,\cdot,R_{i})_{i}$  such that the first order theory $Th(G,\cdot,R_{i})_{i}$ is stable. 

A complete theory $T$ in a language $L$ is said to be {\em stable} if there do {\em not} exist an $L$-formula $\delta(x,y)$, a model $M$ of $T$ and $a_{i},b_{i}\in M$ for $i<\omega$ such that $M\models \delta(a_{i},b_{j})$ iff $i<j$. There are various equivalent conditions to this definition, for example involving counting types. But checking stability of a given theory may not be easy. 
Stable groups may often arise as groups definable in  models $M$ of  stable theories $T$ (equipped with some or all of the structure induced from $M$). Typical examples are where $T$ is the theory of algebraically closed fields of a given characteristic, and we get the class of algebraic groups, or where $T$ is the theory of differentially closed fields of characteristic zero, and we get the class of differential algebraic groups (in the sense of Kolchin). 

As mentioned earlier stable groups support a very nice theory of ``genericity" for definable sets. If $X$ is a definable subset of $G$ then we say $X$ is {\em left-generic} in $G$ if finitely many left translates of $X$ by elements of $G$ cover $G$. Likewise for right generic. $G$ is said to be {\em connected} if $G$ has no proper definable subgroup of finite index. Note that connectedness of $G$ passes to any elementarily equivalent group. 
\begin{Fact} Assume $G$ to be a stable group, and $X,Y$ definable subsets of $G$
Then
\newline
(i) $X$ is left generic iff $X$ is right generic.
\newline
(ii) If $X\cup Y$ is generic, then one of $X$, $Y$ is generic. 
\newline
(iii) $G$ is connected if any only if there is no definable subset $X$ of $G$ such that both $X$ and $G\setminus X$ are generic.
\newline
(iv) For any formula $\phi(x,y)$ of the language of $G$ there is $n$ such that for any $b\in G$, $\phi(x,b)$ (or rather the set it defines) is generic if and only if some $n$ left translates of $\phi(x,b)$ cover $G$ iff some $n$ right translates of $\phi(x,b)$ cover $G$.
\newline
(v) Assume $G$ to be saturated and let $G_{0}$ be a small elementary substructure of $G$. Then $X$ is generic in $G$ iff for every $g\in G$, $(g\cdot X) \cap G_{0} \neq \emptyset$ iff for every $g\in G$, $(X\cdot g)  \cap G_{0}\neq \emptyset$.
\end{Fact}

\noindent
It follows from Sela's Theorem (A) that Fact 1.1 holds for a free nonabelian group. It would be interesting to know if any of (i)-(v) can be proved directly before knowing stability of $T_{fg}$.

Lying behind all of this is the theory of forking in stable structures: if $M$ is a model of a stable theory $T$ and $a\in M$, and $C\subseteq B \subseteq M$ we have the notion: ``$a$ is independent from $B$ over $C$", or ``$tp(a/B)$ does not fork over $C$". The definition is as follows. Write $tp(a/B)$ as $p_{B}(x)$. Then $p_{B}(x)$ does not fork over $A$ if whenenever $\{B_{i}:i<\omega\}$ is indiscernible over $C$ with $B_{0} = B$ then $\{p_{B_{i}}(x):i<\omega\}$ is consistent. (Assuming $M$  saturated and $B,C$ small.) Passing again to a stable group $G$, if $g\in G$ and $A$ is a set of parameters from $G$ then we say that $g$ is generic over $A$, or $tp(g/A)$ is generic, if every formula in $tp(g/A)$ (namely every $A$-definable set containing $g$) is generic in $G$. We have the following forking-theoretic characterizations of genericity:
\begin{Fact} Let $G$ be a saturated stable group.
\newline
(i) $g\in G$, $A\subset G$ small. Then $g$ is generic over $A$ if and only if whenever $h\in G$, and $g$ is independent from $h$ over $A$, then $g$ is independent from $h\cdot g$ over $A$. 
\newline
(ii) Let $p(x)$ be a complete $1$-type over $G$. Then $p(x)$ is generic iff for each $g\in G$, $g\cdot p$ does not fork over $\emptyset$ iff for each $g\in G$, $p\cdot g$ does not fork over $\emptyset$.
\end{Fact}

If the stable group $G$ is connected, then (by Fact 1.1(iii)) there is over any set of parameters a unique complete generic type of an element of $G$. In particular there is a unique generic type, say $p_{0}$ over $\emptyset$ and the generic types over sets of parameters are just the nonforking extensions of $p_{0}$. 

This theory holds with obvious modifications for stable  transitive group actions (stable homogeneous spaces), in particular for stable principal homogeneous spaces, and it will be used below. We repeat a useful fact:
\begin{Fact} Let $G$ be a stable group. Let $H$ be a subgroup of $G$ defined over $A$. Let $g\in G$ be generic over $A$. Then $g/H$ is generic in $G/H$ over $A$, and $g$ is generic in $g\cdot H$ over $A$ together with the canonical parameter of $g\cdot H$.
\end{Fact}

There are several other notions which are important for this paper, such as $T^{eq}$, strong types, type-definable groups, canonical bases etc. , but we refer the reader to \cite{Pillay-book} for example.

\vspace{2mm}
\noindent
I would like to thank Zlil Sela for several discussions, and  to Bestvina and Feign for
allowing me to mention their results. I have benefited from comments from 
several model-theorists, but special thanks are due to Gregory Cherlin for pointing out 
some mistakes and asking some pertinent questions when I spoke on this 
topic at the University of  Lyon I in June 2006.

\section{Genericity and definability in the free group}

Recall our notation: $T_{fg}$ denotes the (complete)
theory of free noncommutative groups in the language of
groups. $F$ denotes a free group $F_{n}$ for some $n\geq 2$,
namely a {\em standard model} of $T_{fg}$. $G$ will denote
an arbitrary, possibly saturated model of $T_{fg}$.

We will typically let $e_{1},...,e_{n}$ denote
free generators of $F_{n}$. By a word in $F_{n}$ we mean a
finite sequence of ``bits" $e_{i}$ and
$e_{i}^{-1}$ for $i=1,..,n$ (which represents of course the
product of these elements). A reduced word is a word in
which for no $i$ is  $e_{i}$ is next to $e_{i}^{-1}$.  So
every member of $F$ is represented by a unique reduced word
(the identity being represented by the empty word). By
``word" we will usually mean reduced word unless we say
otherwise. If $w$ is a word, then by an embedded subword of
$w$ we mean a word $w'$ together with an embedding of $w'$
in $w$ as a sequence of consecutive elements. If we do not
specify the embedding we simply say subword.

As a matter of notation if $A$ is a subset of a group $H$ we let $<A>$ denote the subgroup
of $H$ generated by $A$.

The following definition is due to Bestvina and Feighn.
\begin{Definition} A subset $X$ of $F$ is {\em negligible}  if there 
is a natural number $N$ such that for every $\epsilon > 0$, there is 
a cofinite subset $X'$ of $X$, such that for each $w\in X'$, there 
are $N$ pairs $w_{1}, w_{1}'$, $w_{2}, w_{2}'$,...,$w_{N}, w_{N}'$ of 
proper embedded subwords of $w$, such that
\newline
(i) for each $i$, $w_{i}' = w_{i}$ or $w_{i}^{-1}$, as words.
\newline
(ii) for each $i$, $w_{i} \neq w_{i}'$ as embedded subwords.
\newline
(iii) The $w_{i}$'s and $w_{i}'$'s cover all but $\epsilon$ of $w$. 
That is, the number of elements (bits) of $w$ which are not in any of 
the embedded subwords $w_{i}$ or $w_{i}'$ is $\leq \epsilon\cdot 
length(w)$.
\end{Definition}

\begin{Remark} (i) The nonnegligible subsets of $F$ form a proper 
ideal in the Boolean algebra of subsets of $F$. Namely the union of 
two negligible sets is negligible, a subset of a negligible set is 
negligible, and $F$ itself is nonnegligible.
\newline
(ii) If $X\subseteq F$ is negligible then so is any left or right 
translate of $X$ by an element of $F$.
\newline
(iii) Any cyclic subgroup of $F$ is negligible.
\newline
(iv) The commutator subgroup $[F,F]$, and its complement, are both negligible.
\end{Remark}
{\em Proof.}  Clear. But in (ii) for example, when passing to a 
translate $g\cdot X$, and considering a word $w\in X$, there may be 
some cancellation when we pass to the reduced word $gw$, which must 
be taken account of.

\vspace{5mm}
\noindent
The following substantial result has been announced by Bestvina and Feighn \cite{B-F}:
\begin{Proposition} Let $X$ be a definable subset of $F$. Then either 
$X$ or $F\setminus X$ is negligible ( and by Remark 2.2 (i), not both).
\end{Proposition}

We now give several consequences.

\begin{Proposition} (i) Let $X$ be a definable subset of
$F$. Then $X$ is non-negligible if and only if $X$ is
generic.
\newline
(ii) For any formula $\phi(x,y)$ of $L$, there is some
formula $\psi(y)$ of $L$ such that for any (free group) $F$
and $b\in F$, $\phi(x,b)(F)$ is non-negligible if and only
if $F\models \psi(b)$.
\end{Proposition}
{\em Proof.}   (i) The right to left direction is by
Remark 2.2 (i) and (ii).
\newline
Left-to-right. Suppose $X$ is nongeneric.
Then by Fact 1.1 (iv), $F\setminus X$ is generic,
hence non-negligible by the right to left direction, hence
$X$ is negligible by Proposition 2.3.
\newline
(ii) follows from (i) by Fact 1.1 (iv). 

\begin{Proposition} The free group has a unique
generic type, and in particular is connected. That is, for
any model $G$ of $T_{fg}$, (a) for any definable subset $X$
of $G$, precisely one of $X$, $G\setminus X$ is generic, and
(b) $G$ has no proper definable subgroup of finite index.

\end{Proposition}
{\em Proof.} Clear from Fact 1.1.

\vspace{5mm}
\noindent
Let $p_{0}(x)\in S_{1}(T_{fg})$ be the unique generic
type of $T_{fg}$ (over $\emptyset$). 

\vspace{5mm}
\noindent
Actually Proposition 2.5 can be also deduced just from
Theorems A and B of Sela, using the following elementary observation
of Bruno Poizat from around 25 years ago. (This result appeared in an early draft of
of \cite{Poizat} but was for some reason omitted in the published version.)
\begin{Lemma} Let $X$ be a definable subset of $F_{\omega} =
<e_{n}:n=1,2,....>$. Suppose $X$ is generic. Then for all
but finitely many $n$, $e_{n}\in X$.
\end{Lemma}
{\em Proof.} If $g_{1}X \cup .... \cup g_{s}X = G$, let
$r$ be such that the parameters in the formula defining $X$
as well as $g_{1},..,g_{r}$ are words in $e_{1},..,e_{r}$
and their inverses.  Let $i>r$. So $e_{i}\in g_{t}X$ for some $t$,
whence $g_{t}^{-1}e_{i}\in X$. But there is an automorphism
of $F_{\omega}$ fixing each of $e_{1},..,e_{r}$ and taking
$g_{t}e_{i}$ to $e_{i}$. So $e_{i}\in X$.

\vspace{5mm}
\noindent
Lemma 2.6 implies that there are no two disjoint definable
generic subsets of $F_{\omega}$. By Theorem (B), $F_{\omega}$ is a model of $T_{fg}$,
which by Theorem (A) is stable. So again we conclude from Fact 1.1 the  uniqueness of the generic type, and connectedness.

\vspace{5mm}
\noindent
In fact Lemma 2.6 gives a bit more:
\begin{Corollary}  (i) In $F_{\omega} = <e_{i}:i<\omega>$,
the sequence $(e_{i}:i<\omega)$ is a Morley sequence in the generic type 
$p_{0}$.
\newline
(ii) In $F = F_{n}$, $(e_{1},..,e_{n})$ is an
independent set of realizations of $p_{0}$.
\end{Corollary}
{\em Proof.} (i) Let $X$ be a generic definable subset of
$F_{\omega}$,  defined over $e_{1},...,e_{r}$ say. We have
seen that $e_{i}\in X$ for some $i>r$. But then (by
automorphism) $e_{r+1}\in X$. This shows that
$tp(e_{r+1}/e_{1},..,e_{r})$ is generic (i.e. =
$p_{0}|\{e_{1},..,e_{r}\}$) which is what we wanted to prove.
\newline
(ii) follows as Theorem (B) implies that $F_{n}$ is an
elementary substructure of $F_{\omega}$ (under the canonical
embedding).

\vspace{5mm}
\noindent
Note that for any $m>1$, the set of $mth$ powers in $F$ is a
negligible definable set, hence non generic (hence non
generic in any model). (Alternatively by Corollary 2.7
we see that the generic type $p_{0}(x)$ implies ``$x$ is not an mth power"
for all $m$.)
However the $mth$ power map is
injective, so the free group could not be superstable (as
pointed out by Poizat in the early draft of \cite{Poizat}). Noting also that every
element is the product of a square and a cube, we have:
\begin{Proposition} (i) Let $G$ be a saturated model of
$T_{fg}$ and $A$ a small subset of $G$. Then $\{g\in G:
tp(g/A)$ is not generic\} is not a subgroup of $G$.
\newline
(ii) The generic type $p_{0}$ of $G$ does not have weight
$1$. That is, it is NOT the case that for any set $A$ of
parameters, forking on realizations of $p_{0}|A$ is an
equivalence relation.
\end{Proposition}
{\em Proof.} (ii) follows from (i) by standard manipulations.

\vspace{5mm}
\noindent
We now pass to definable subgroups. Before we begin note
that the free group is centreless.
\begin{Lemma} (i) Any nontrivial abelian subgroup of a (nonabelian) free
group $F$ is cyclic.
\newline
(ii) Any nontrivial negligible subgroup of a (nonabelian) free group $F$
is abelian, so cyclic.
\end{Lemma}
{\em Proof.} (i) is obvious. (ii) is also reasonably
obvious. For example if $H$ is a nonabelian subgroup of $F$
then $H$ contains a free group $F'$. But clearly $F'$ is
nonnegligible, hence so is $H$.

\begin{Corollary} Let $F$ be a (nonabelian) free group.
\newline
(i) Any proper
definable subgroup of  $F$ is abelian.
\newline
(ii) The following three classes of subgroups of $F$
coincide:
\newline
(a) maximal abelian subgroups,
\newline
(b) $\{C(a): a\neq 1, a\in F\}$,
\newline
(c) The maximal proper definable subgroups.
\newline
  Moreover for $a\neq 1$, $C(a)$ is the unique maximal
abelian subgroup of $F$ containing $a$.
\newline
(iii) For any $a,b\in F$ different from $1$ either $C(a) =
C(b)$, or $C(a)\cap C(b) = \{1\}$.
\newline
(iv) The above (i), (ii) and (iii) hold in any model $G$ of
$T_{fg}$.
\newline
(v) If $G$ is a model of $T_{fg}$ and $a\neq 1$ is in $G$
then $C(a)$ is self-normalizing.
\end{Corollary}
{\em Proof.} (i) As $F$ is connected any proper definable
subgroup of $F$ is of infinite index, hence nongeneric,
hence nonnegligible (by 2.3). Now use Lemma 2.9.
\newline
(ii) Let $B$ be a maximal abelian subgroup of $F$. Let $a\in
B$, then $C(a)$ is a proper definable subgroup of $F$
containing $B$ but by (i) is abelian, hence coincides with
$B$.
\newline
Let $a\neq 1$. Then $C(a)$ is proper, definable, and by (i)
abelian. Moreover any proper definable subgroup containing
$C(a)$ is abelian (by (i)) hence coincides with $C(a)$.
\newline
Let $B$ be maximal proper definable. So $B$ is abelian (by
(i)). Any abelian subgroup containing $B$ is definable, so
coincides with $B$.
\newline
The moreover clause is contained in the proof above.
\newline
(iii) Suppose $c\in C(a)\cap C(b)$, $c\neq 1$. Then $C(c)$
is maximal abelian and contains $a$ and $b$. By (ii) $C(a) =
C(b) = C(c)$.
\newline
(iv) This follows by transfer.
\newline
(iv) It is enough to work in a standard model $F$. Then
$C(a)$ is cyclic, with generator $u$ say. It is then easy to
finf $v\in F$ such that $u^{v}$ is not a power of $u$. So
$C(a)$ is not normal. So the normalizer of $C(a)$ is a
proper definable subgroup of $F$, hence by (ii) coincides
with $C(a)$.

\begin{Corollary} The free group is definably simple, namely
has no proper nontrivial definable normal subgroups.
\end{Corollary}

A {\em simple bad} group of finite Morley rank is a
simple group
$G$ of finite Morley rank (definable in some ambient
structure) such that the Borels of $G$, namely maximal
connected solvable subgroups of $G$, are nilpotent. If $B$
is such a Borel, then it is known that $B$ is
self-normalizing and that distinct Borels are disjoint over
$\{1\}$. Hence the free group resembles such a simple bad
group of finite Morley rank, where we interpret 
``Borel" in a free group as a maximal abelian subgroup. But note that in  free groups, 
these ``Borels" are not connected (and this will introduce an interesting twist to the proof of
Proposition 3.2 in the next section).

In a simple bad group of finite
Morley rank, any easy calculation shows that the union of
the conjugates of a Borel $B$ is generic in $G$. However this
fais in the free group:

\begin{Lemma} Let $G$ be a model of $T_{fg}$ and $B = C(a)$
for
$a\neq 1$. Then $\cup_{g\in G}B^{g}$ is not generic in $G$.
\end{Lemma}
{\em Proof.} It suffices to work in a standard model $F$. We
will just consider for simplicity the case $B = C(e_{1})$
where
$e_{1}$ is one of the generators of $F$. In this case
clearly $B = <e_{1}>$. Then any nonidentity element of
$\cup_{g\in F}B^{g}$ is {\em when put in reduced form} of
the form
$w^{-1}e^{\pm m}w$ for some $w$ (possibly empty) and $m\geq
1$. Moreover, for any $k<\omega$ there are clearly only
finitely many such reduced words of length $k$ with $m = 1$.
So taking $N = 2$ we see that $\cup_{g}B^{g}$ is negligible,
hence nongeneric.

\vspace{5mm}
\noindent
{\bf Exercise.} Show, more generally that for any nongeneric definable subset $X$ of $G\models T_{fg}$, $\cup_{g\in G}X^{g}$ is nongeneric.

\section{Geometric stability}
We begin in the context of an arbitrary complete stable
theory $T$, working in a saturated model ${\bar M}$ of $T$.
In fact we work freely in ${\bar M}^{eq}$. Our aim is to prove that the (theory of the) free group is not $CM$-trivial.  As motivation we
first discuss $1$-basedness (or modularity).
\begin{Definition} $T$ is not $1$-based (or $T$ is
$1$-ample), if either of the following equivalent conditions
hold:
\newline
(i) there are tuples $a,b$ such that $tp(a/acl(b))$ forks
over $acl(a)\cap acl(b)$.
\newline
(ii) after possibly adding parameters, there are $a,b$ such
that $a$ forks with $b$ over $\emptyset$ but $acl(a)\cap
acl(b) = acl(\emptyset)$,
\newline
(iii) There are $a,B$ such that $Cb(stp(a/B))$ is not
contained in $acl(a)$.
\end{Definition}

It is a basic fact that if $G$ is a group definable in a
$1$-based stable theory then $G$ is abelian-by-finite. Hence
the (theory of the) free group is not $1$-based. The existence of non
$1$-based theories of finite Morley rank (or $U$-rank) in
which no infinite field (or even group) is definable is a
nontrivial fact. However it is easy to find such non
$1$-based structures if we drop the finite Morley rank
condition. One such is the {\em free pseudoplane} which we
discuss briefly as the technology is related to what we do
with the free group.

The free pseudoplane is the theory with one binary relation
$I$ axiomatized by
\newline
(i) $I$ is symmetric and irreflexive,
\newline
(ii) for all $x$ there are infinitely many $y$ such that
$I(x,y)$, and
\newline
(iii) there are no $I$-loops of length $\geq 3$, namely ther
do not exist distinct $x_{0},x_{1},..,x_{n}$ with $n\geq 2$
such that $I(x_{i},x_{i+1})$ for $i<n$ and $I(x_{n},x_{0})$.

The free pseudoplane  is a complete $\omega$-stable
theory. The unique $1$-type over $\emptyset$ has Morley rank
$\omega$. If $M$ is a model and $a\neq b\in M$ then the
Morley rank of $tp(b/a)$ is the length of the shortest
$I$-path from $a$ to $b$ if there is one, or $\omega$
otherwise. Moreover all types over parameters are
stationary. These are all easy to verify. We claim that
the non $1$-basedness of the free pseudoplane is witnessed
in the following strong form:
\newline
{\em Claim.} Let $a,b$ be such that $I(a,b)$. Then $a$ forks
with $b$ over $\emptyset$ but $acl(a)\cap acl(b) = acl(\emptyset)$ where $acl(-)$ is computed in
$M^{eq}$.
\newline
{\em Proof.} The forking is clear as $RM(tp(b)) = \omega$
but $RM(tp(b/a)) = 1$.
\newline
Now suppose for a contradiction that there is $e\in
acl^{eq}(a)\cap acl^{eq}(b)\setminus acl^{eq}$. So $b$ forks
with $e$, witnessed by a formula $\phi(x,e)$ where we may
assume that $\phi(x,e)$ implies $e\in acl(x)$. So
$\phi(x,e)$ has Morley rank $N$. It follows that
\newline
(*) for all
$b',b''$ realising
$\phi(x,e)$ the shortest path joining $b'$ and $b''$ is at
most $N$.
\newline
However, let $a_{0}, b_{0}, a_{1}, a_{2},...$ be
chosen as follows: $a_{0} = a$, $b_{0} = b$, $a_{i+1}\neq
a_{i}$ has same strong type as $a_{i}$ over $b_{i}$, and
$b_{i+1}\neq b_{i}$ has the same strong type as $b_{i}$ over
$a_{i+1}$. As $e\in acl(a)\cap acl(b)$, each $b_{i}$
realizes $\phi(x,e)$. But as there are no loops, the
shortest path between $b_{0}$ and $b_{n}$ has length $2n$,
giving a contradiction.

\vspace{5mm}
\noindent
The nonabelianness of the free group gives a canonical
configuration witnessing non-basedness (as in Definition
3.1). But there is another configuration witnessing non
$1$-basedness which has  the same character as in the Claim above
for the free pseudoplane. We discuss this now. Putting
the two pseudoplanes together will give non $CM$-triviality as we explain
subsequently.

Let us take $F$ to be a (finitely generated) free group on at
least
$4$ generators
$\{e_{1}, e_{2}, e_{3}, e_{4},...\}$. We can consider $F$ as
an elementary substructure of a saturated model $G$, but in
fact we will work in the standard model $F$.
We will need the following lemma, which is left as an exercise:
\begin{Lemma} Fix $k$ and let $Y\subset F$ be the set of
words of the the form
$e_{1}^{k}(e_{4}^{-1}e_{1}^{k})(e_{4}^{-2}e_{1}^{k})...(e_{4}^{-(n-1)}e_{1}^{k})e_{4}^{n(n-1)/2}$,
as $n$ varies.
Then $Y$ is nonneglible.
\end{Lemma}

One of our main results is:

\begin{Proposition}  Let $B = C(e_{1}) = <e_{1}>$. Let
$l$ be a canonical parameter for the translate $e_{2}\cdot
B^{e_{3}}$ of $B^{e_{3}}$. Work over $e_{1}$ (namely add a
constant for
$e_{1}$). THEN $acl^{eq}(e_{2}) \cap acl^{eq}(l) =
acl(\emptyset)$, but $e_{2}$ forks with $l$ over $\emptyset$.
\end{Proposition}
{\em Proof.} As in the statement of the Proposition we
work over $e_{1}$.
The forking is clear as
$e_{2}$ is generic over
$\emptyset$ but is in the nongeneric definable set
$e_{2}B^{e_{3}}$ which has canonical parameter $l$.

For the rest: Assume for a contradiction that
$d\in acl^{eq}(e_{2})\cap acl(l)\setminus acl(\emptyset)$.
(Note that a priori we know nothing about imaginaries in
$F^{eq}$.) So $e_{2}$ forks with $d$ and this is witnessed by
a nongeneric formula $\phi(x,d)$ satisfied by $e_{2}$ and
without loss implying that $d\in acl^{eq}(x)$. The coset
$e_{2}B^{e_{3}}$ is a principal homogeneous space for $B^{e_{3}}$.

\begin{Lemma} (i)  $tp(e_{2}/l)$ is a generic type of
$e_{2}B^{e_{3}}$, and 
\newline
(ii) $\phi(x,d)\wedge x\in l$ defines, up to a
nongeneric subset of $e_{2}B^{e_{3}}$, a union of orbits
under the $kth$ powers of $B^{e_{3}}$ for some $k$.
\end{Lemma}
{\em Proof.} (i) follows from Fact 4.3.
\newline
(ii) Let $r = stp(e_{2}/l) = tp(e_{2}/acl^{eq}(l))$,  a stationary generic type of $e_{2}B^{e_{3}}$. One knows that $r(x)$ is determined by the data ``$x$ is generic over $l$ in $e_{2}B^{e^{3}}$" together with the orbit of $e_{2}$ under the connected component of $B^{e_{3}}$. However the connected component of $B^{e_{3}}$ (in a saturated model) is simply the intersection of the $kth$ powers of $B^{e_{3}}$ for all $k$ (as $B^{e_{3}}$ is torsion-free abelian). Now $\phi(x,d)\in tp(e_{2}/acl^{eq}(l))$, so it follows from the above comments together with compactness that for some $k$, the set of elements of $e_{2}B^{e_{3}}$ satisfying $\phi(x,d)$ is, up to a nongeneric set, a union of orbits under the $kth$ powers of $B^{e_{3}}$. 

\vspace{2mm}
\noindent
Fix $k$ as in Lemma 3.4(ii). The technical core of this paper is contained in the following.
\begin{Lemma} Let $c_{0} = e_{2}$, $g_{0} = e_{3}$ and
$l_{0} = l$. Let $c_{i}, g_{i}, l_{i}$ for $i>0$ be defined
by
\newline
(a) $g_{i} = e_{3}\cdot e_{4}^{i(i+1)/2}$, and 
\newline
(b) $c_{i} = c_{i-1}\cdot(e_{1}^{k})^{g_{i-1}}$,
\newline
(c) $l_{i}$ is the canonical parameter for $c_{i}\cdot
B^{g_{i}}$.
\newline
Then for all $i\geq 0$.
\newline
(i) $c_{i}$ and $g_{i}$ are independent generic over $e_{4}$, and
$tp(c_{i}/l_{i})$ is a generic of $c_{i}B^{g_{i}}$.
\newline
(ii) $c_{i}$ satisfies $\phi(x,d)$,
\newline
(iii) $tp(l_{i},d) = tp(l_{0},d)$
\newline
(iv) $c_{i+1}$ and $g_{i}$ are independent generic over $e_{4}$.
\newline
(v) $tp(c_{i+1}/l_{i})$ is a generic of $c_{i}B^{g_{i}}$
and $c_{i+1}$ satisfies $\phi(x,d)$.
\end{Lemma}
{\em Proof.} Let us first do the case $i=0$. (i) and (ii) and (iii) are already given to us
(using the Claim).
\newline
(iv) We have (i) for $i=0$, so $c_{0}$ is generic over $\{g_{0},e_{4}\}$, hence $c_{0}$ is generic over
$\{(e_{1}^{k})^{g_{0}}, g_{0}, e_{4}\}$, hence $c_{1} = c_{0}(e_{1}^{k})^{g_{0}}$ is generic over $\{(e_{1}^{k})^{g_{0}}, g_{0}, e_{4}\}$. In particular $c_{1}$ is generic over $\{g_{0},e_{4}\}$ giving (iv).
\newline
(v) As $c_{1}$ is in the same orbit as $c_{0}$ under the
$kth$ powers of $B^{g_{0}}$ we see that
$\models\phi(c_{1},d)$. It is easy to see that $c_{0}$ is
independent from $(e_{1}^{k})^{g_{0}}$ over $l_{0}$, hence
$c_{1}$ is a generic of $l_{0}$.
\newline
Now assume (i) - (iv) are true for $i$ and we will prove them for $i+1$.
\newline
(i) The induction assumption gives that $c_{i+1}$ and $g_{i}$ are generic, independent over $e_{4}$. But $g_{i+1} = g_{i}e_{4}^{i+1}$, hence  $c_{i+1}$ and $g_{i+1}$ are generic, independent over $c_{4}$. So by 4.3, $tp(c_{i+1}/l_{i+1})$ is a generic of $c_{i+1}B^{g_{i+1}}$.  
\newline
(ii) is given by the induction hypothesis.
\newline
(iii) By (iv) of the induction hypothesis, the fact that $d\in acl(c_{i+1})$, part  (i), and stationarity of the generic type of $F$ we see that $tp(g_{i},c_{i+1}/d) = tp(g_{i+1},c_{i+1}/d)$ and hence $tp(l_{i}/d) = tp(l_{i+1}/d)$. 
\newline
(iv) and (v) are proved as in the case $i=0$. 
\newline
Lemma 3.5 is proved.

\vspace{2mm}
\noindent
Note that (after cancelling) $c_{n} = e_{2}e_{3}^{-1}(e_{4}^{-1}e_{1}^{k})(e_{4}^{-2}e_{1}^{k}).....(e_{4}^{-(n-1)}e_{1}^{k})e_{4}^{n(n-1)/2} e_{3}$.
But each $c_{n}$ realizes  $\phi(x,d)$ which is nongeneric. Hence $\{c_{n}:n<\omega\}$ is negligible. This is clearly a contradiction with Lemma 3.2. 

\begin{Remark}  Note that from Proposition 3.3 we obtain that if $G$ is a saturated model of $T_{fg}$ and $e_{1},c,g\in G$ are any independent generics, $B = C(e_{1})$ and $l$ is a canonical parameter for $b\cdot T^{g}$, then after naming $e_{1}$, $acl(c)\cap acl(l) = acl(\emptyset)$. (Simply because $tp(e_{1},e_{2},e_{3})$ in $F$ equals $tp(e_{1},c,g)$ in $G$.)
\end{Remark}

\vspace{5mm}
\noindent
We now give the notion of $CM$-triviality, working again freely in ${\bar M}^{eq}$ for ${\bar M}$ a saturated model of $T$.
\begin{Definition}  $T$ is not $CM$-trivial (or $T$ is $2$-ample) if there are tuples $c,b,a$ such that
\newline
(i) $acl(a,b)\cap acl(a,c) = acl(a)$.
\newline
(ii) $a = Cb(stp(c/a))$ and $b = Cb(stp(c/ab))$.
\newline
(iii) $a\notin acl(b)$. 
\end{Definition}

Equivalent statements, which do not mention canonical bases are:
\newline
(I). $T$ is not $CM$-trivial if, possibly after adding parameters, there are $a,b,c$ such that $acl(a)\cap acl(b) = acl(\emptyset)$, $acl(a,b)\cap acl(a,c) = acl(a)$, $c$ is independent from $b$ over $a$, and $a$ forks with $c$ over $\emptyset$.

\vspace{2mm}
\noindent
(II) $T$ is not $CM$ trivial if there exist $a,b,c$ such that $a$ is independent from $b$ over $c$, but $a$ forks with $b$ over $acl(a,b)\cap acl(c)$

\vspace{5mm}
\noindent
The notion was introduced by Hrushovski \cite{Hrushovski} where he also proved the equivalence of the three versions for strongly minimal theories. The proof goes through for arbitrary stable theories. None of the definitions is particularly memorable, but version (I) is stated in a manner that suggests natural strengthenings. In fact that is what we did in \cite{Pillay2}, introducing the notion $n$-ample for any $n\geq 1$. As pointed out by several people, including David Evans and Ikuo Yoneda our definition needed some additional fine tuning:
\begin{Definition}  Let $n\geq 1$. Then $T$ is $n$-ample if (after possibly naming parameters) there are $a_{0},..,a_{n}$ such that
\newline
(i) $acl(a_{0})\cap acl(a_{1}) = acl(\emptyset)$ and $acl(a_{0},a_{1},...a_{i-1},a_{i})\cap acl(a_{0},a_{1},..,a_{i_1},a_{i+1}) = acl(a_{0},..,a_{i-1})$ for $1\leq i <n$, and 
\newline
(ii) $a_{i+1}$ is independent from $\{a_{0},...,a_{i-1},a_{i}\}$ over $a_{i}$ for all $1\leq i < n$, and 
\newline
(iii) $a_{0}$ forks with $a_{n}$ over $\emptyset$.
\end{Definition}

A stable field is $n$-ample for all $n$ (\cite{Pillay2}). In \cite{Evans} David Evans found, for each $n$ a stable $n$-ample theory which is moreover a reduct of a trivial $1$-based theory (so interprets no infinite groups).  We believe the free group to be non $3$-ample. 

On the other hand our proof \cite{Pillay1} of non $CM$-triviality (or $2$-ampleness) of bad groups of finite Morley rank readily generalizes to the free  group, using Proposition 3.3:
\begin{Proposition} $T_{fg}$ is non $CM$-trivial.
\end{Proposition}
{\em Proof.} We will be brief. Fix a saturated model $G$ of $T_{fg}$.  Fix $e_{1}\in G$ generic, and let $T = C(e_{1})$. Add a constant for $e_{1}$ (namely work over $e_{1}$). Let $a,g,b,c\in G$ be independent generics. Let $G_{a} = \{(h,h^{a}):h\in G\}$, and let $P$ be a canonical parameter for the coset $(b,c)\cdot G_{a}$.

Let $(T^{g})_{a} = \{(h,h^{a}):h\in T^{g}\}$ (a subgroup of $G_{a}$), and let $l$ be a canonical parameter for the coset $(b,c)\cdot (T^{g})_{a}$. 

We want to check that the triple $(P,l,(b,c)$ witnesses non $CM$-triviality of $T_{fg}$ as in Definition 3.7, namely
\newline
(i) $acl(P,l) \cap acl(P,(b,c)) = acl(P)$,
\newline
(ii) $P = Cb(stp((b,c)/P)$ and $l = Cb(stp((b,c)/P,l))$.
\newline
(iii) $P\notin acl(l)$.

\vspace{2mm}
\noindent
Proof of (i). $(b,c)\cdot G_{a}$ is a $P$-definable PHS for $G_{a}$, hence also for $G$ (as $G$ is isomorphic to $G_{a}$ via $h\to (h,h^{a})$). So fixing a point $d\in P$ generic over the data $\{b,c,g\}$ gives a $\{P,d\}$-definable bijection between $(b,c)\cdot G_{a}$ and $G$, which takes $(b,c)$ to $b'$ say and $(b,c)(T^{g})_{a}$ to $b'\cdot T^{g}$. Moreover $b',g$ are generic independent over $P$. Letting $l'$ be a canonical parameter for $b'\cdot T^{g}$ we see from Proposition 3.3 that $acl(b')\cap acl(l) = acl(\emptyset)$. This implies easily that $acl(P,(b,c))\cap acl(P,l) = acl(P)$ as required.

\vspace{2mm}
\noindent
(ii) is routine, as $(b,c)$ is a generic point of $(b,c)\cdot G_{a}$ over $P$, and also a generic point of $(b,c)\cdot (T^{g})_{a}$ (which has canonical parameter $l$) over $\{P,l\}$.

\vspace{2mm}
\noindent
(iii). For $a'\in G$, let $P' = (b,c)\cdot G_{a'}$ and $l' = (b,c)\cdot (T^{g})_{a'}$ (as elements of $G^{eq}$). If $a'\in a\cdot T$, $a'\neq a$ then $l' = l$ but $P'\neq P$. We can choose (infinitely many) such $a'$ such that $tp(a',g,b,c) = tp(a,g,b,c)$, so $tp(P',l) = tp(P,l)$ for infinitely many distinct $P'$.

\section{Questions and problems}
We list some further problems (with some commentaries), some of which may be  settled by the literature, or even be obvious after a little reflection. 

\begin{Problem} Suppose $G$ is a model of $T_{fg}$. Is the free product of $G$ and $\Z$ an elementary extension of $G$?
\end{Problem} 

\noindent
{\em Comment.}  By definability of the generic type $p_{0}$, Problem 4.1 is equivalent to: Let $F = F_{n}$ be some/any nonabelian free group. Let $\phi({\bar z},x,y)$ be a formula (in the language of groups). Then there are terms (or words), $t_{1}({\bar z},y),..,t_{r}({\bar z},y)$ such that for any $\bar m$ from $F$, if $\exists x \phi({\bar m},e_{n+1},x)$ holds in $F_{n+1}$, then there is $i\leq r$ such that
$F_{n+1}\models \phi({\bar m},e_{n+1},t_{i}({\bar m},e_{n+1}))$.
\newline
I would imagine this result to be contained in the proof of relative quantifier elimination for $T_{fg}$. 
\newline
Note that (a positive answer to) problem 4.1 implies that the free product of any model $G$ of $T$ with any any free group is an elementary extension of $G$.

\vspace{2mm}
\noindent
A related question is:
\begin{Problem} Let $G$ be  a model of $T_{fg}$ and $G_{1}, G_{2}$ two elementary extensions of $G$. Is the free product of $G_{1}$ and $G_{2}$ over $G$ an elementary extension of each of $G_{1}, G_{2}$?
\end{Problem}

\begin{Problem} Let $B$ be a ``Borel" in a model $G$ of $T_{fg}$ (namely $B$ is a maximal proper definable subgroup). Does $B$ have $U$-rank $1$. More generally, are the subsets of $B^{n}$ definable in $G$ precisely those definable in the structrure $(B,\cdot)$.
\end{Problem}

\noindent
{\em Comment.} This should be known. It is probably enough to prove this for {\em standard} $B$, namely $B$ defined over some $F = F_{n}$  ($n>1$). In this case one may try to carry out Ehrenfeucht-Fraisse games in $F$ in its relational language to obtain the desired conclusion.

\begin{Problem} Describe the $U$-rank $1$ types in $(T_{fg})^{eq}$, and more generally the superstable (type)-definable sets.
\end{Problem}

\begin{Problem} Does Proposition 2.4 (equivalently Bestvina-Feighn's Proposition 2.3) hold for $F_{\omega}$ under the same definition of negligible?
\end{Problem}

\noindent
{\em Comment.} Note that under Definition 2.1 applied to $F_{\omega} = <e_{i}:i<\omega>$, the set $(e_{i}:i<\omega)$ would be non-negligible.
\newline
Problem 4.5 is related to:

\begin{Problem} Let $\phi(x)$ be a formula over some nonabelian free group $F = F_{n}$ which is nongeneric (defines a nongeneric set in $F$). So for every $m\geq n$, $\phi(F_{m})$ is negligible. Is there an $N$ as in Definition 2.1 which works for all $\phi(F_{m})$. 
\end{Problem}

\begin{Problem} Prove that no infinite field is interpretable in $T_{fg}$.
\end{Problem}

\begin{Problem} Prove that $T_{fg}$ is not $3$-ample.
\end{Problem}

\begin{Problem} Describe the saturated models of $T_{fg}$.
\end{Problem}

\noindent
{\em Comment.} For example the saturated models of $Th(\Z, +)$ are of the form  ${\hat \Z} \oplus\Q^{\kappa}$ (where ${\hat\Z}$ is the profinite completion, or pure-injective hull, of $\Z$). So in particular a $\kappa$-saturated model of $T_{fg}$ will contain a free product of $\kappa$-copies of ${\hat \Z} \oplus\Q^{\kappa}$

\end{document}